\documentclass{amsart}
\usepackage{subcaption}
\usepackage{graphicx}
\usepackage[vflt]{floatflt}
\usepackage{rotating}
\usepackage{pdflscape}
\usepackage[dvipsnames]{xcolor}
\usepackage{colortbl}
\usepackage{mathtools}
\usepackage{amssymb}
\usepackage{blkarray}
\usepackage{tikz-cd}
\usepackage[all]{xy}
\usepackage{multicol}
\usepackage{enumitem}
\usepackage{setspace}
\usepackage{todonotes}
\usepackage{comment}
\usepackage{etoolbox}
\usepackage[hypertexnames=false,
    colorlinks,
    citecolor=red,
    linkcolor=blue,
    urlcolor=red,
    bookmarksopen,
    backref=page]{hyperref}
\usepackage[depth=3]{bookmark}
\usepackage{charter}
\usepackage{newtxmath}

\numberwithin{equation}{section}

\newtheorem{theorem}{Theorem}[section]
\newtheorem{lemma}[theorem]{Lemma}
\newtheorem{proposition}[theorem]{Proposition}

\newtheorem{corollary}[theorem]{Corollary}
\newtheorem*{thm*}{Theorem}

\theoremstyle{definition}
\newtheorem{example}[theorem]{Example}

\newtheorem{definition}[theorem]{Definition}

\newtheorem{mainthm}{Theorem}

\newcommand{\be}{\begin{equation}}
\newcommand{\ee}{\end{equation}}
\newcommand{\bes}{\begin{equation*}}
\newcommand{\ees}{\end{equation*}}

\newcommand{\cA}{\mathcal{A}}

\newcommand{\cD}{\mathcal{D}}

\newcommand{\cH}{\mathcal{H}}

\newcommand{\cK}{\mathcal{K}}

\newcommand{\cQ}{\mathcal{Q}}
\newcommand{\cR}{\mathcal{R}}

\newcommand{\cX}{\mathcal{X}}

\newcommand{\cZ}{\mathcal{Z}}

\newcommand{\bB}{\mathbb{B}}
\newcommand{\bC}{\mathbb{C}}
\newcommand{\bD}{\mathbb{D}}

\newcommand{\bN}{\mathbb{N}}

\newcommand{\bR}{\mathbb{R}}
\newcommand{\bS}{\mathbb{S}}
\newcommand{\bT}{\mathbb{T}}

\newcommand{\lip}{\langle}
\newcommand{\rip}{\rangle}
\newcommand{\ip}[1]{\lip #1 \rip}
\newcommand{\bip}[1]{\big\lip #1 \big\rip}

\newcommand{\re}{\operatorname{Re}}
\newcommand{\im}{\operatorname{Im}}

\newcommand{\sD}{\scr{D}}

\newcommand{\qand}{\quad\text{and}\quad}

\newcommand{\ba}{{\mathbf{a}}}

\def\cK{\mathcal{K}}

\def\scr{\mathscr}

\def\bi{\begin{itemize}}
\def\ei{\end{itemize}}

\def\bpm{\begin{pmatrix}}
\def\epm{\end{pmatrix}}
\newcommand{\bsm}{\left[\begin{smallmatrix}}
\newcommand{\esm}{\end{smallmatrix} \right] }

\def\ba{\begin{eqnarray*}}
\def\ea{\end{eqnarray*}}

\IfFileExists{accents.sty}{\usepackage{accents}}{}

\hypersetup{
    colorlinks=true,
    linkcolor=blue,
    citecolor=red,
    urlcolor=blue
}
\usepackage[nameinlink,noabbrev]{cleveref}
\crefname{section}{Section}{Sections}
\crefname{lemma}{Lemma}{Lemmas}
\crefname{theorem}{Theorem}{Theorems}

\newcommand{\rank}{\operatorname{rank}}
\newcommand{\Hol}{\operatorname{Hol}}

\begin{document}

\title{Complete characterization of cyclic polynomials on polyballs}
\author{Jeet Sampat}
\makeatletter
\patchcmd{\@setauthors}{\endtrivlist}{%
  \par\smallskip
  {\normalfont\footnotesize
   Department of Mathematics, University of Manitoba\par
   Email: \href{mailto:Jeet.Sampat@umanitoba.ca}
        {\texttt{Jeet.Sampat@umanitoba.ca}}}%
  \endtrivlist
}{}{}
\makeatother

\begin{abstract}
    	We give a complete characterization of cyclic polynomials in the total-degree Dirichlet-type spaces on polyballs $\Omega$, that is, Cartesian products of  Euclidean unit balls in arbitrary finite dimension. For a stable polynomial $P$ with a nonempty distinguished boundary zero set, the critical index is $D-\kappa$, where $D$ is the complex dimension of the polyball and $\kappa$ is the largest exponent for which the distinguished boundary zero-set of $P$ supports a probability measure assigning mass $O(t^\kappa)$ to every Kor\'anyi ball of radius $\sqrt t$. In this case, $P$ is cyclic in the Dirichlet-type space $\sD_s(\Omega)$ with parameter $s \in \bR$ exactly when $s\le D-\kappa$. On the other hand, if $P$ is stable and the distinguished boundary zero-set is empty then $P$ is cyclic in $\sD_s(\Omega)$ for all $s \in \bR$.
    	
    	Our proof combines a new Hardy kernel construction with an induction argument on the sum of block degrees of $P$, supplied by stability-preserving polar derivatives. A sharp estimate for the surface measure of neighborhoods of the zero set closes the induction at the critical index $D - \kappa$, while the probability measure yields a bounded functional	annihilating every polynomial multiple of $P$ above it. This gives a unified proof for the unit ball and the polydisk, where the critical index depends only on the real dimension of the boundary zero set. A direct corollary of our main result is that $P$ is cyclic in the Drury--Arveson space exactly when it is stable and its boundary zero-set is empty or has real dimension at most two.
\end{abstract}

\nobreakdash
\maketitle

\section{Introduction}

\subsection{Cyclicity}

Let $T$ be a bounded linear operator on a separable Banach space $\cX$. A vector $v \in \cX$ is said to be \emph{$T$-cyclic} if the collection of vectors $T^jv$, $j=0,1,\dots$, spans a dense subspace of $\cX$. A common difficulty in the characterization of cyclic vectors is the need to develop techniques specific to the space $\cX$ and sometimes even to the operator $T$ under consideration. A classical illustration of this phenomenon is the Smirnov--Beurling theorem \cite{Beu49,Smi32}, which characterizes all cyclic vectors of the canonical shift on the complex sequence space $\ell^2(\bN)$ by unitarily identifying it with the shift operator $f \mapsto zf$ on the Hardy space $H^2(\bD)$ of analytic power series on the unit disk with square-summable coefficients. The Smirnov--Beurling theorem then states that a function in $H^2(\bD)$ is shift-cyclic if and only if it is \emph{outer}, giving a function theoretic characterization of this operator-theoretic property. This cornerstone of function theory continues to motivate cyclicity research in other analytic function spaces. See the author's survey \cite{Sam25} for more background on cyclicity. In this article, we introduce a new Hardy kernel method that provides a complete characterization of cyclic polynomials on polyballs in every dimension, unifying the unit ball and the polydisk. In particular, we answer a conjecture of Vavitsas and Zarvalis \cite{VZ24} in the positive with the unit ball characterization, which in turn provides a characterization of cyclic polynomials in the Drury--Arveson space.

\subsection{Polyball domains}

Fix $d_1,\ldots,d_r \in \bN$ and write $D = \sum_k d_k$. Let $\bB_{d_k}$ be the Euclidean unit ball in $\bC^{d_k}$ for each $1\le k\le r$, and let $\bS_{d_k} := \partial \bB_{d_k}$ be the complex unit sphere. The corresponding \emph{polyball domain} and \emph{polysphere} are, respectively,
\begin{equation*}
    \Omega=\prod_{k=1}^r\mathbb B_{d_k} \qand
    \Sigma=\prod_{k=1}^r\mathbb S_{d_k}.
\end{equation*}
Moreover, the set $\Sigma$ is the distinguished boundary of $\Omega$, and we write $\sigma$ for the product of normalized surface measures on it. Clearly, the unit ball $\bB_D$ and the unit polydisk $\bD^D$ are both polyballs. We invoke the standard multi-index notations:
\begin{equation*}
    \bN_0 = \bN \cup \{0\}, \quad \alpha^{(k)} = (\alpha^{(k)}_1, \dots, \alpha^{(k)}_{d_k}) \in \bN_0^{d_k}, \quad \alpha=(\alpha^{(1)},\ldots,\alpha^{(r)}) \in \bN_0^D,
\end{equation*}
\begin{equation*}
    \alpha! = \prod_{k = 1}^r \prod_{j = 1}^{d_k}\alpha^{(k)}_j!, \quad z^\alpha = \prod_{k = 1}^r \prod_{j = 1}^{d_k} \big({z^{(k)}_j}\big)^{\alpha^{(k)}_j}, \qand |\alpha| = \sum_{k = 1}^r \sum_{j = 1}^{d_k} \alpha^{(k)}_j.
\end{equation*}
Consider the family of \emph{Dirichlet-type spaces} $\sD_s(\Omega)$, $s\in\mathbb R$, consisting of all functions $f(z)=\sum_\alpha f_\alpha z^\alpha\in\Hol(\Omega)$ such that
\begin{equation*}
    \|f\|_s^2= \sum_{\alpha \in \bN_0^D}(1+|\alpha|)^s \prod_{k=1}^r \frac{(d_k-1)!\,\alpha^{(k)}!}{(d_k-1+|\alpha^{(k)}|)!} |f_\alpha|^2<\infty.
\end{equation*}
Each $\sD_s(\Omega)$ is a Hilbert space under the inner-product associated with $\left\|\cdot\right\|_s$, and the set of polynomials is densely contained in each $\sD_s(\Omega)$.

The \emph{Hardy space} on $\Omega$ is the space $H^2(\Omega) := \sD_0(\Omega)$, and it is the central object that our cyclicity results will be based around. For a single ball, $\sD_s(\bB_D)$ is the usual Dirichlet-type scale: $s=0$ gives $H^2(\Omega)$, $s=D-1$ gives the Drury--Arveson space $\cH^2_D$ up to an equivalent norm, and $s=-1$ gives a norm equivalent to the unweighted Bergman norm. Similarly, for the polydisk, we recover the \emph{total-degree} Dirichlet-type spaces as considered in \cite{NZ26} for the bidisk.

A polynomial $P$ is said to be \emph{cyclic} if its invariant subspace
\begin{equation*}
    [P]_s=\overline{\{Pf \ : \ f\in\sD_s(\Omega)\}}^{\|\cdot\|_s}
\end{equation*}
is equal to $\sD_s(\Omega)$. $P$ is said to be \emph{stable} if it has no zeros in $\Omega$. For any $P$, we write $E(P)=\{\eta\in\Sigma:P(\eta)=0\}$ for its distinguished boundary zero-set.

Fix a polyball domain $\Omega \subset \bC^D$ and introduce the metric
\begin{equation*}
    \varrho(\eta, \xi) := \max_{1 \leq k \leq r} |1 - \ip{\eta^{(k)}, \xi^{(k)}} |^{\frac{1}{2}}, \quad \eta,\xi \in \Sigma.
\end{equation*}
In fact, the triangle inequality for $\varrho$ extends to all of $\overline \Omega$: for $x, y, v \in \overline \bB_d$,
\begin{equation*}
    1 - \ip{x,y} = [(1 - \ip{x,v}) + (1 - \ip{v,y})] - [\ip{x-v,y-v} + 1 - |v|^2].
\end{equation*}
Then, use the fact that $|x-v|^2 + 1 - |v|^2 \leq 2 \re (1 - \ip{x,v})$ along with Cauchy--Schwarz. Several of our estimations revolve around the \emph{Kor\'anyi balls} $$B_t(\eta) := \{ \xi \in \Sigma \ : \ \varrho(\eta, \xi) < \sqrt{t}\} \subset \Sigma$$ and the \emph{tube} $U_t(E):=\bigcup_{\eta\in E}B_t(\eta)$, where $E$ is a nonempty compact subset of $\Sigma$. Let $$\varrho(\omega,E)=\min_{\eta\in E}\varrho(\omega,\eta)$$ so that $U_t(E)=\{\omega\in\Sigma:\varrho(\omega,E)<\sqrt t\}$.

\subsection{Main results and outline of proofs}

We use the standard notation $A \lesssim B$ to denote that $A \leq cB$ for some constant $c > 0$ independent of the varying parameters, and write $A \asymp B$ if both $A \lesssim B$ and $B \lesssim A$ hold.

\begin{definition}
Let $E(P)\ne\varnothing$ and $\kappa\ge0$. We say that $P$ satisfies the \emph{$\kappa$-measure property} if there exists a Borel probability measure $\mu$ supported on $E(P)$ such that
        \begin{equation}\label{eq:polynomial.measure}
            \mu\left(B_t(\eta)\right) \lesssim t^\kappa, \quad \forall 0 < t \leq 1/2, \, \eta \in \Sigma. \tag{$\mu_\kappa$}
        \end{equation}
\end{definition}

Let $\kappa_*(P)$ be the supremum of the exponents in \eqref{eq:polynomial.measure}. If $E(P)=\varnothing$, we use the convention $\dim_\bR E(P)=\kappa_*(P)=-\infty$. Lemma~\ref{lem:cap.cover.measure} shows that the supremum is attained and its value is intrinsic to $E(P)$. However, using \eqref{eq:intrinsic.exponent.formula}, one can compute it via some smooth parametrization of $E(P)$. Combining this with Lemma \ref{lem:product.trace} immediately yields our first main result.

\begin{mainthm}\label{mainthm:covering-measure.estimates}
Suppose $P$ is a nonzero polynomial on $\Omega$ and $E(P)\ne\varnothing$. Then $\kappa_*(P)\in\tfrac12\bN_0$, $\kappa_*(P)<D$, and \eqref{eq:polynomial.measure} holds at $\kappa=\kappa_*(P)$. Moreover, we have the tube estimate
    \begin{equation*}
        \sigma(U_t(E(P)))\asymp t^{D-\kappa_*(P)},\quad \forall0<t\le1/2,
    \end{equation*}
    Any measure $\mu$ as in \eqref{eq:polynomial.measure} with $\kappa=\kappa_*(P)$ satisfies the Hardy-trace property
    \begin{equation}\label{eq:Hardy-trace}
        \int_{E(P)} \left\lvert f((1-t) \eta) \right\rvert^2 \, d\mu(\eta) \lesssim t^{\kappa_*(P) - D} \|f\|^2_0, \quad \forall 0 < t \leq 1/2, \, f \in H^2(\Omega).
    \end{equation}
\end{mainthm}

The upper tube estimate supplies the approximation at $D - \kappa_*(P)$ to derive cyclicity of $P$ in $\sD_s(\Omega)$, whereas the Hardy-trace estimate gives a bounded functional annihilating $[P]_s$ above this index. Thus, we arrive at our characterization result.

\begin{mainthm}\label{mainthm:characterization.of.cyclicity}
    $P \not\equiv 0$ is cyclic in $\sD_s(\Omega)$ if and only if it is stable and $s \leq D - \kappa_*(P)$.
\end{mainthm}

The proof of cyclicity is through an induction argument on the block-total-degree of $P$, i.e., $\sum_k \deg_k P$, and each $\deg_k P$ is the total-degree of $P$ with variables from the $\bB_{d_k}$ block. Lemma~\ref{lem:polars} and Proposition~\ref{prop:product.approximation} show that the \emph{polar derivatives}
\begin{equation*}
    Q_k := \begin{cases}
        (\deg_k P) P - \sum_{j = 1}^{d_k} z^{(k)}_j \frac{\partial P}{\partial z^{(k)}_j}, \quad &\forall k \ : \ \deg_kP>0,\\
        0, \quad &\forall k \ : \ \deg_k P = 0,
    \end{cases} 
\end{equation*}
have smaller block-total-degree and satisfy $E(Q_k) \subseteq E(P)$ whenever $\deg_kP>0$. Repeated polar differentiation in active blocks eventually
leads to a nonzero constant, the base case of the induction. For $Q=\prod_{k\,:\,\deg_kP\,>\,0}Q_k$, using the upper tube estimate in Theorem \ref{mainthm:covering-measure.estimates}, we show that $Q^{D+1} \in [P]_s$ whenever $s \leq D - \kappa_*(P)$. Since $Q$ is cyclic if and only if $Q^{D+1}$ is cyclic, this helps us complete the induction step.

Recall that a nonempty real algebraic set $E$ has real dimension $m$ if it is the largest value for which some projection onto $m$ real coordinates has image with nonempty interior in $\bR^m$. In the classical domains $\bB_D$ and $\bD^D$ the exponent $\kappa_*(P)$ for any stable $P$ is completely determined by its real dimension $m_P := \dim_\bR(E(P))$ (see Lemma~\ref{lem:classical.covering}). Thus, the following corollaries are immediate from Theorem \ref{mainthm:characterization.of.cyclicity}.

\begin{corollary}\label{cor:ball.characterization}
    If $P$ is stable on $\bB_D$, then $\kappa_*(P) = m_P/2$ and thus
    \begin{equation*}
        P \text{ is cyclic in } \sD_s(\bB_D) \quad\iff\quad m_P \leq 2(D - s).
    \end{equation*}
    In particular, $P$ is cyclic in the Drury--Arveson space $\cH^2_D$ if and only if $m_P \leq 2$.
\end{corollary}

\begin{corollary}\label{cor:polydisk.characterization}
    If $P$ is stable on $\bD^D$, then $\kappa_*(P) = m_P$ and thus
    \begin{equation*}
        P \text{ is cyclic in } \sD_s(\bD^D) \quad\iff\quad m_P \leq D - s.
    \end{equation*}
\end{corollary}

The one-variable endpoint is classical~\cite{BS84}. In dimension
$D=2$, Kosi\'nski and Vavitsas settled the ball case~\cite{KV23}, and
Nailwal and Zalar completed the bidisk case~\cite{NZ26}, following
\cite{TZ26}. For $\bB_D$ in arbitrary dimension, Vavitsas and Zarvalis
proved the necessary condition in
Theorem~\ref{mainthm:characterization.of.cyclicity} for stable $P$
with $m_P>0$~\cite{VZ24}. We recover it directly from the Hardy-trace
estimate supplied by~\eqref{eq:polynomial.measure} and prove the
matching sufficiency at the critical index. In particular, this
gives the characterization in $\cH^2_D$ for $D\geq 3$. For related results on the ball see \cite{Vav23,VWZ26}, and for the polydisk see \cite{BCLSS15II,BKKLSS16,Berg18}.

On the unit ball, Aleman, Perfekt, Richter, Sundberg and Sunkes \cite[Theorem~3.6]{APRSS24} showed that cyclicity passes from a stable polynomial $P$ to another stable
polynomial $Q$ whenever $E(Q) \subseteq E(P)$. Our method achieves this concretely by starting with $P$ directly: stability-preserving polar derivatives give lower-degree polynomials $Q$ and a modified \L{}ojasiewicz estimate, then Hardy kernels and analytic measure estimates close the induction at the critical index.

\bookmarksetupnext{level=1}
\subsubsection*{\texorpdfstring{\textbf{Disclosure of AI usage}}{Disclosure of AI usage}}

ChatGPT-6 assisted in preliminary brainstorming and with manuscript preparation. The author developed the polar reduction and the kernel estimate method independently of this assistance. Initial tests were done for $\cH^2_3$ which led to a very technical version of Proposition \ref{prop:product.approximation} that completed the induction argument. The author subsequently simplified and generalized the method to work for polyballs in any dimension. The author takes full responsibility for the content of this article and has no affiliation with OpenAI.

\section{Properties of the Dirichlet-type spaces}\label{sec:properties.Dirichlet-type.space}

For simplicity, fix a polyball $\Omega \subset \bC^D$ and write $\sD_s = \sD_s(\Omega)$ for each $s \in \bR$. Let
\begin{equation*}
 h_\alpha=\prod_{k=1}^r
       \frac{(d_k-1)!\,\alpha^{(k)}!}
            {(d_k-1+|\alpha^{(k)}|)!},
 \quad w_s(\alpha)=(1+|\alpha|)^sh_\alpha.
\end{equation*}
The monomial integral in each block of $\Sigma$ gives $\int_\Sigma z^\alpha\overline{z^\beta}\,d\sigma(z) = \delta_{\alpha\beta}h_\alpha.$ Consequently the norm at $s=0$ is the Hardy norm, and the
homogeneous terms $f_n=\sum_{|\alpha|=n}f_\alpha z^\alpha$
are orthogonal in every $\sD_s$. In particular,
\begin{equation}\label{eq:homogeneous.norm}
 \|f\|_s^2=\sum_{n=0}^\infty (n+1)^s\|f_n\|_0^2,
 \quad
 \|f\|_0^2=\sup_{0<\rho<1}
       \int_\Sigma |f(\rho\zeta)|^2\,d\sigma(\zeta).
\end{equation}

\begin{enumerate}[wide,itemsep=2mm]
\item \textbf{Reproducing kernels.}
The normalized monomials $z^\alpha/\sqrt{w_s(\alpha)}$ are orthonormal to each other. For $\rho=\max_k|u^{(k)}|<1$, $u\in\Omega$, the multinomial identities
give
\begin{equation*}
 \sum_{\alpha \in \bN_0^D}\frac{|u^\alpha|^2}{w_s(\alpha)}
 \le\sum_{n=0}^\infty (n+1)^{-s}
       \binom{n+D-1}{D-1}\rho^{2n}<\infty.
\end{equation*}
The Cauchy--Schwarz inequality therefore bounds point evaluation at each $u \in \Omega$. The same estimate is uniform on compact subsets, so coefficient partial sums converge locally uniformly as well as in norm. This identifies $\sD_s$ with a Hilbert space of holomorphic functions that is a reproducing kernel Hilbert space, and proves that the set of polynomials is densely contained in each $\sD_s$.
At $s=0$, summing the binomial series yields the
Hardy--Szeg\H{o} kernel
\begin{equation}\label{eqn:Hardy--Szego.kernel}
 S(z,u)=\prod_{k=1}^r
       \big(1-\langle z^{(k)},u^{(k)}\rangle\big)^{-d_k},
 \quad \|S(\cdot,u)\|_0^2=S(u,u).
\end{equation}

\item \textbf{Shifts and multipliers.} Let $e_{k,j}$ denote the multi-index of $z_j^{(k)}$ and note
\begin{equation*}
 \frac{w_s(\alpha+e_{k,j})}{w_s(\alpha)}
   =\left(\frac{|\alpha|+2}{|\alpha|+1}\right)^s
      \frac{\alpha_j^{(k)}+1}{|\alpha^{(k)}|+d_k}
   \le2^{\max(s,0)}.
\end{equation*}
Thus each coordinate shift $S_{k,j} : f \mapsto z^{(k)}_j f$ is bounded, and every polynomial is a multiplier. Recall that a function $H$ is a multiplier of $\sD_s$ if the multiplication operator $M_H : f \mapsto Hf$ is bounded on $\sD_s$ into itself.

\item \textbf{Weak convergence.} Suppose $f_k \in \sD_s$ and $f \in \Hol(\Omega)$ are such that
\begin{equation*}
    \lim_{k \to \infty} f_k(u) = f(u); \quad \forall u \in \Omega, \qand \sup_{k \in \bN} \|f_k\|_s < \infty.
\end{equation*}
By the reproducing kernel property, the span of the evaluation kernels $K_{s,u}$, $u\in\Omega$, is dense in $\sD_s$. Thus, if $\varphi$ is a finite linear sum of these kernel functions, the pointwise convergence of the $f_k$'s shows that $\ip{f_k, \varphi}$ is a convergent sequence for all such $\varphi$. The uniform norm bound on the $f_k$'s extends this convergence to $\ip{f_k, H}$ for each $H \in \sD_s$. In particular, $f_k$ converges weakly, and the Riesz representation theorem identifies the resulting weak limit with $f$. Thus $f\in\sD_s$ and $f_k \rightharpoonup f$ (weak convergence) as $k \to \infty$. Every norm-closed linear subspace is weakly closed, so a closed subspace containing all $f_k$ also contains $f$.

\item \textbf{Cyclicity and polynomial products.}
If $P(u)=0$ at any $u \in \Omega$, then the bounded evaluation at $u$
annihilates $[P]_s$ and equals $1$ on the constant
function $1$. Thus a cyclic polynomial must be stable. Since polynomial multiplication is bounded and polynomials are
dense, $[P]_s$ is also the closure of the polynomial multiples
of $P$. In particular, $Pf\in[P]_s$ for every $f\in\sD_s$. If $[P]_s$
contains a cyclic vector $g$, then it contains the closure of all
polynomial multiples of $g$, and so it equals $\sD_s$. Thus, $P$ is cyclic if and only if $1 \in [P]_s$.

For polynomials $P_1,P_2$, one has $[P_1P_2]_s\subset[P_1]_s\cap[P_2]_s$, so cyclicity of the product implies cyclicity of both factors. Conversely, if both factors are cyclic, choose polynomials $Q_k$ with $P_2Q_k\to1$ in $\sD_s$. Bounded multiplication by $P_1$ gives $P_1P_2Q_k\to P_1$, so $P_1\in[P_1P_2]_s$. This subspace contains the cyclic vector $P_1$ and is the entire space. The same conclusion holds for any finite product $P_1 \cdots P_k$. In particular, for every power $P^k$, $k \in \bN$, we have $P^k$ is cyclic if and only if $P$ is cyclic.

Lastly, we have $s' < s \implies \|f\|_{s'} \leq \|f\|_s$. Thus if $P$ is cyclic at $s$ then the same sequence of polynomials $Q_k$ with $Q_k P \to 1$ in $\sD_s$ shows that $P$ is cyclic in $\sD_{s'}$.

\item \textbf{Radial derivatives and interpolation.}
Let $\mathcal R=\sum_{k=1}^r\mathcal R_k$, where $$\cR_k := \sum_{j = 1}^{d_k} z_j^{(k)} \tfrac{\partial}{\partial z^{(k)}_j}$$ is the radial derivative on $\bB_{d_k}$. The identity $\mathcal R f_n=nf_n$, $n \in \bN$, gives
\begin{equation}\label{eq:radial.norm.identity}
 \|f\|_{2c}=\|(I+\mathcal R)^cf\|_0
      \asymp |f(0)|+\|\mathcal R^cf\|_0, \quad \forall c \in \bN.
\end{equation}
The equivalence follows from $(n+1)^{2c}\asymp n^{2c}$, $n \geq 1$.
For $t \ge 0$ and $0\le\theta\le1$, coefficient-wise
H\"older's inequality gives
\begin{equation}\label{eq:scale.interpolation}
 \|f\|_{\theta t}^2\le\|f\|_0^{2(1-\theta)}\|f\|_t^{2\theta}.
\end{equation}
In particular, if $H\in\Hol(\overline\Omega)$ then $H$ has bounded radial derivatives of every order. By \eqref{eq:radial.norm.identity} and monotonicity of the norms, $H\in\sD_s$ for every $s\in\bR$.

\end{enumerate}

\section{Theorem A: The critical exponent \texorpdfstring{$\kappa_*$}{kappa*}}\label{sec:proof.Theorem.A}

Let $P$ be a nonzero polynomial on a polyball $\Omega$ with $E(P) \ne \varnothing$. The equations $\re P=\im P= 0$ and $1-|z^{(k)}|^2=0$ for all $ 1 \leq k \leq r$ express $E(P)$ as a compact real algebraic set. Applying Nowak's theorem on covering with special cubes \cite[p.~67]{Now09} to $|P|^2$ on $\Sigma$ near $E(P)$ (cf.\ \cite[\S~4.5]{Gro87}), we obtain the chart parametrization
\[
    E(P)=\bigcup_{j=1}^J F_j(\cQ_j),\quad \cQ_j=(-1,1)^{m_j}.
\]
Each $F_j$ is a real-analytic diffeomorphism onto its image and extends real analytically into $\Sigma$ on a neighborhood of the closure of a larger cube $\cQ_j^+\supset\overline{\cQ_j}$. If $m_j = 0$ for any $j$, then we treat $\cQ_j$ as just a single point at the origin. If $r$ is the number of blocks defining $\Omega$, then each $m_j<2D-r$. Indeed, $m_j = 2D-r$ implies $F_j(\cQ_j)$ contains an open subset of $\Sigma$, so $P|_\Sigma\equiv0$ by connectedness and real analyticity, and hence $P\equiv0$ by the maximum modulus principle.
For $x\in\cQ_j^+$ and $y\in\bR^{m_j}$, define
\[
    A_j(x)y=\bigl(\im\bip{DF_j^{(k)}(x)y,F_j^{(k)}(x)}\bigr)_{k=1}^r .
\]
Differentiating $|F_j^{(k)}(x+ty)|^2=1$ at $t = 0$ shows that $\re\bip{DF_j^{(k)}(x)y,F_j^{(k)}(x)}=0$. Since the second derivatives of $F_j$ are bounded on $\cQ_j^+$, Taylor's formula gives
\begin{equation}\label{eq:cap.chart.expansion}
    \left|
        \varrho(F_j(x'),F_j(x))^2 - \|A_j(x)(x'-x)\|_{\max} \right| \lesssim |x'-x|^2,
\end{equation}
where $x,x'\in\cQ_j^+$ and $\|\cdot\|_{\max}$ is the maximum norm on $\bR^r$.

We are now ready to introduce the exponent $\kappa_*$ from the introduction:
\begin{equation}\label{eq:intrinsic.exponent.formula}
    h_j := \max_{x \in \cQ_j} \big(\rank A_j(x) \big) \qand \kappa_* := \frac{1}{2} \max_{1 \leq j \leq J} (m_j + h_j) \in \frac{1}{2} \bN_0.
\end{equation}
As $m_j < 2D - r$ and $h_j \leq r$ for all $1 \leq j \leq J$, it follows that $\kappa_* < D$. Even though $\kappa_*$ is defined through a chart, the following result shows that it is intrinsic to $E(P)$. Consequently, it is intrinsic to $P$ and $\Omega$ together.

\begin{lemma}\label{lem:cap.cover.measure}
The exponent $\kappa_*$ is the largest $\kappa\ge0$ for which a probability measure on $E(P)$ satisfies \eqref{eq:polynomial.measure} with exponent $\kappa$. Moreover, for $U_t(E(P)) = \bigcup_{\eta \in E(P)} B_t(\eta)$,
\begin{equation*}
    \sigma(U_t(E(P)))\asymp t^{D-\kappa_*},\quad 0<t\le1/2.
\end{equation*}
\end{lemma}

\begin{proof}
By the standard Kor\'anyi ball estimate (see, e.g., \cite[Lemma~4.6]{Zhu05}) applied to each block $1 \leq k \leq r$, it is easy to see that $\sigma(B_t(\eta))\asymp t^D$, $0<t\le 1/2$. For any finite positive measure $\lambda$ on $\Sigma$, symmetry and Tonelli's theorem give
\begin{equation}\label{eq:ball.mass.identity}
    \int_\Sigma\lambda(B_t(\zeta))\,d\sigma(\zeta)
    =\int_\Sigma\sigma(B_t(\eta))\,d\lambda(\eta)
    \asymp t^D\lambda(\Sigma).
\end{equation}
Let $\lambda_j$ be the Lebesgue measure on $\bR^{m_j}$ (with point mass one if $m_j=0$). Let $\nu$ be the pushforward measure on $\Sigma$, given by
\[
    \nu(B)=\sum_{1 \leq j \leq J}\lambda_j\bigl(\cQ_j^+\cap F_j^{-1}(B)\bigr), \quad \forall \text{ Borel } B \subseteq \Sigma.
\]
At $x\in\cQ_j$ write $y=y_0+y_1 \in \ker A_j(x) \oplus \ker A_j(x)^\perp$, and put $q=\rank A_j(x)$. If $|y_0|<c\sqrt t$ and $|y_1|<ct$, then $x+y\in\cQ_j^+$ for sufficiently small fixed $c$ determined by the distance from $\overline \cQ_j$ to $\bR^{m_j} \setminus \cQ_j^+$. Using \eqref{eq:cap.chart.expansion} and boundedness of $A_j$, we get
\[
    \varrho(F_j(x+y),F_j(x))^2
    \lesssim |y_1|+|y|^2\lesssim(c+c^2)t.
\]

Choosing a smaller $c$, if needed, so that the R.H.S. above is $< t$, we have
    \begin{equation*}
        \cA_j := \{ y_0 + y_1 \ : \ |y_0| < c \sqrt{t}, \, |y_1| < ct \} \subseteq \cQ_j^+ - x, \quad F_j(x + \cA_j) \subseteq B_t(F_j(x)).
    \end{equation*}
Note that $\lambda_j(\cA_j) \asymp t^{(m_j+q)/2}$ as $\dim \ker A_j(x) = m_j - q$ and $\dim \ker A_j(x)^\perp = q$. Since $(m_j + q)/ 2 \leq \kappa_*$ by definition and $0 < t < 1$, we get that $t^{\kappa_*} \lesssim \lambda_j(\cA_j)$. Using the definition of $\nu$ and the above containment for $\cA_j$ we get $t^{\kappa_*} \lesssim \nu(B_t(F_j(x)))$. For $\zeta\in U_t(E(P))$, some $B_t(F_j(x))$ lies in $B_{4t}(\zeta)$ by the triangle inequality on $\varrho$, and so
\begin{equation*}
    t^{\kappa_*} \lesssim \nu(B_{4t}(\zeta)), \quad \forall \zeta \in U_t(E(P)).
\end{equation*}
Hence \eqref{eq:ball.mass.identity} and the standard Kor\'anyi ball estimate give
\[
    t^{\kappa_*}\sigma(U_t(E(P)))\lesssim\int_\Sigma\nu(B_{4t}(\zeta))\,d\sigma(\zeta) = \int_\Sigma \sigma(B_{4t}(\zeta)) d\nu(\zeta)\lesssim t^D.
\]
Note that the Kor\'anyi ball estimate in the last estimate above only applies to $t \leq 1/8$. If $t > 1/8$, then use the fact that $\sigma(\Sigma) = 1$ and so the bound is simply $\nu(\Sigma)\leq (8^D \nu(\Sigma)) t^D$. Thus, we have shown that
\begin{equation}\label{eq:sigma.tube.bound}
    \sigma(U_t(E(P))) \lesssim t^{D - \kappa_*}
\end{equation}

Let $j$ be such that $m_j+h_j=2\kappa_*$ and pick a small closed cube $\cK\subset\cQ_j$ about a point $x$ with $\rank A_j(x) = h_j$. For $y = y_0 + y_1$ as before, $\|A_j(x)y\|_{\max}\gtrsim|y_1|$ on $\cK$. This follows from the continuity of the $h_j$-th singular value if $h_j>0$, and from $y_1=0$ if $h_j=0$. Suppose $x,x+y\in\cK$ and $F_j(x+y)\in B_t(F_j(x))$ both hold. Let $\eta = F_j(x)$ and $\xi = F_j(x+y)$. Then, by the choice of $\eta, \xi$ we have
\begin{equation*}
    |\xi^{(k)}- \eta^{(k)}|^2 = 2 \re (1 - \ip{\xi^{(k)},\eta^{(k)}}) \leq 2 |1 - \ip{\xi^{(k)}, \eta^{(k)}}| < 2t, \quad \forall 1 \leq k \leq r
\end{equation*}
The squared Euclidean distance in $\bC^D$ is the sum of the L.H.S. above, and so $$|F_j(x+y) - F_j(x)|^2 \leq 2rt.$$ Note that $F_j$ is bi-Lipschitz on $\cK$ by the inverse function theorem applied to the $m_j$ independent real coordinates, and so this means that $|y| \lesssim \sqrt{t}$. Since $F_j(x+y) \in B_t(F_j(x))$, put $x' = x+y$ in \eqref{eq:cap.chart.expansion} to note that $\|A_j(x)y\|_{\max}\lesssim t$. Thus $|y_0|\lesssim\sqrt t$ and $|y_1|\lesssim t$. For the probability measure
\begin{equation}\label{eq:def.mu}
    \mu(B)=\frac{\lambda_j(\cK\cap F_j^{-1}(B))}{\lambda_j(\cK)}, \quad \forall \text{ Borel } B \subseteq \Sigma,
\end{equation}
the volume of the set of such $y_0 + y_1$ gives $$\mu(B_t(F_j(x)))\lesssim t^{(m_j+h_j)/2}=t^{\kappa_*}.$$ Now, any ball $B_t(\eta)$ meeting $F_j(\cK)$ is contained in $B_{4t}(F_j(x))$ for a point $F_j(x)$ in their intersection. This gives \eqref{eq:polynomial.measure}, including $m_j=0$ by the convention above.

For any probability measure $\lambda$ supported on $E(P)$ and satisfying \eqref{eq:polynomial.measure} with exponent $\kappa$, we have
\[
    t^D\lesssim t^\kappa\sigma(U_t(E(P)))\lesssim t^{D+\kappa-\kappa_*}.
\]
Indeed, the upper bound follows from \eqref{eq:sigma.tube.bound}. For the lower bound, note that if $\zeta \notin U_t(E(P))$ then $B_t(\zeta) \cap E(P) = \varnothing$, and so $\lambda(B_t(\zeta)) = 0$. Otherwise, by \eqref{eq:polynomial.measure} we get $\lambda(B_t(\zeta)) \lesssim t^\kappa$. Thus, \eqref{eq:ball.mass.identity} gives
\begin{equation*}
    t^D \lesssim \int_{U_t(E(P))} \lambda(B_t(\zeta)) d\sigma(\zeta) \lesssim t^\kappa \sigma(U_t(E(P))),
\end{equation*}
and so $\kappa \le \kappa_*$. Taking $\lambda=\mu$ as in \eqref{eq:def.mu} and $\kappa=\kappa_*$ gives the reverse estimate.
\end{proof}

\begin{lemma}\label{lem:product.trace}
    If $E$ is a nonempty compact subset of $\Sigma$, and $\mu$ is a Borel probability measure supported on $E$ satisfying \eqref{eq:polynomial.measure}, then
    \begin{equation*}
        \int_{E} \left\lvert f((1-t) \eta) \right\rvert^2 \, d\mu(\eta) \lesssim t^{\kappa - D} \|f\|^2_0, \quad \forall 0 < t \leq 1/2, \, f \in H^2(\Omega).
    \end{equation*}
\end{lemma}

\begin{proof}
    Fix $0<t\le 1/2$ and, for $\eta,\omega\in\Sigma$, put
    \begin{equation*}
        p_t(\eta,\omega)=[t(2-t)]^D
        |S((1-t)\eta,\omega)|^2,
    \end{equation*}
    where $S$ is the Hardy--Szeg\H o kernel \eqref{eqn:Hardy--Szego.kernel}. Let $f$ be a polynomial and set $z=(1-t)\eta$. Since
    $fS(\cdot,z)\in H^2(\Omega)$, the reproducing property gives
    \begin{equation*}
        f(z)S(z,z)
        =\langle fS(\cdot,z),S(\cdot,z)\rangle_0
        =\int_\Sigma f(\omega)|S(\omega,z)|^2\,d\sigma(\omega).
    \end{equation*}
    From \eqref{eqn:Hardy--Szego.kernel}, we have
    $S(z,z)=[t(2-t)]^{-D}$. Dividing the above equation by
    $S(z,z)$ and using $|S(\omega,z)|=|S(z,\omega)|$ gives
    \begin{equation}\label{eq:f_1-t.reproducing}
        f((1-t)\eta)
        =\int_\Sigma f(\omega)p_t(\eta,\omega)\,d\sigma(\omega).
    \end{equation}
    Taking $f=1$ shows that $p_t$ has integral $1$ in the second
    variable. The kernel formula gives
    $$S((1-t)\eta,\omega)=\overline{S((1-t)\omega,\eta)},$$
    so $p_t$ is symmetric and has integral $1$ in the first
    variable as well.

    Suppose $\xi\in B_t(\eta)$. For each $k$, the extended triangle inequality for $\varrho$ on $\overline\Omega$ shows that the two quantities
    \begin{equation*}
        |1-(1-t)\langle\eta^{(k)},\omega^{(k)}\rangle|^{1/2}
        \qand
        |1-(1-t)\langle\xi^{(k)},\omega^{(k)}\rangle|^{1/2}
    \end{equation*}
    differ by at most
    $|1-\langle\eta^{(k)},\xi^{(k)}\rangle|^{1/2}\le\sqrt t$.
    Each is at least $\sqrt t$ because
    $|\langle x,\omega^{(k)}\rangle|\le1$ for
    $x=\eta^{(k)},\xi^{(k)}$. Consequently they are both comparable, and the
    product kernel formula gives
    $p_t(\eta,\omega)\asymp p_t(\xi,\omega)$ uniformly in $\omega$.

    Invoke the Kor\'anyi ball estimate $\sigma(B_t(\eta))\asymp t^D$ once again. Averaging the comparison of $p_t(\eta,\omega)$ and $p_t(\xi, \omega)$ over
    $\xi\in B_t(\eta)$ gives
    \begin{equation*}
        p_t(\eta,\omega)
        \lesssim \frac{1}{\sigma(B_t(\eta))}
                  \int_{B_t(\eta)}p_t(\xi,\omega)\,d\sigma(\xi).
    \end{equation*}
    Since $\xi\in B_t(\eta)$ if and only if $\eta\in B_t(\xi)$,
    Tonelli's theorem yields
    \begin{align*}
        \int_\Sigma p_t(\eta,\omega)\,d\mu(\eta)
        &\lesssim t^{-D}\int_\Sigma\int_{B_t(\eta)}
            p_t(\xi,\omega)\,d\sigma(\xi)\,d\mu(\eta)\\
        &=t^{-D}\int_\Sigma p_t(\xi,\omega)
            \mu(B_t(\xi))\,d\sigma(\xi)\\
        &\lesssim t^{\kappa-D},
    \end{align*}
    where the last inequality follows from
    \eqref{eq:polynomial.measure} and the normalization of $p_t$
    in the first variable. By Cauchy--Schwarz in \eqref{eq:f_1-t.reproducing} and the
    normalization in the second variable,
    \begin{equation*}
        |f((1-t)\eta)|^2
        \le\int_\Sigma |f(\omega)|^2p_t(\eta,\omega)\,d\sigma(\omega).
    \end{equation*}
    Integrating against $\mu$ and using the preceding estimate
    proves \eqref{eq:Hardy-trace} for polynomials. For fixed $t$,
    convergence in the Hardy norm implies uniform convergence
    on $(1-t)\Sigma$, so polynomial density gives the result
    for every $f\in H^2(\Omega)$.
\end{proof}

\begin{proof}[\textbf{Proof of Theorem~\ref{mainthm:covering-measure.estimates}}]
    The conclusions follow directly from both these lemmas.
\end{proof}

\section{Theorem B: Characterization of cyclicity}

Let $P$ be a nonconstant stable polynomial on $\Omega$ and define
\begin{equation*}
    \nabla_k P = \bigl(\partial_1^{(k)}P,\dots,\partial_{d_k}^{(k)}P\bigr) \qand
    \cR_k = z^{(k)} \cdot \nabla_k,
\end{equation*}
where $x \cdot y = \sum_j x_j y_j$. For each block $\bB_{d_k}$, let $n_k=\deg_kP$ and define $$Q_k=n_kP-\mathcal R_kP, \quad Q=\textstyle \prod_{k \, : \, Q_k\not\equiv0}Q_k.$$ A block is called \emph{active} if $n_k>0$ and \emph{inactive} if $n_k=0$. Our induction argument for the proof of Theorem \ref{mainthm:characterization.of.cyclicity} relies on the fact that the \emph{polar derivative} $Q_k$ of $P$ in an active block preserves stability and lowers the degree in that block. The following lemma establishes this property and the accompanying gradient estimates by applying Harris's range-inclusion
theorem on a single ball \cite[Theorem~2]{Harris96} to each block separately. For related polar derivative results, see \cite{deBruijn47,Leake22}.

\begin{lemma}\label{lem:polars}
    If $n_k>0$ then $Q_k$ is stable, $\deg_kQ_k\leq n_k-1$, $\deg_j Q_k \leq n_j$ for $j \neq k$,
    \begin{equation}\label{eq:block.gradient.bound}
        |\nabla_kP|\le |Q_k| \qand n_k|P|\le2|Q_k| \quad\text{on }\overline\Omega.
    \end{equation}
    In particular, $E(Q_k)\subset E(P)$ and $\nabla_kP(\eta)=-Q_k(\eta)\overline{\eta^{(k)}}$, for all $ \eta\in E(P)$.
\end{lemma}

\begin{proof}
    The degree assertions on $Q_k$ follow by cancellation of block degree $n_k$ terms.
    Fix $k$ with $n_k>0$ and arbitrarily chosen $z_j \in \bB_{d_j}$ for each $j \neq k$, and let $$p(x) = P(z_1, \dots, z_{k - 1}, x, z_{k + 1}, \dots, z_r), \quad x \in \bB_{d_k}.$$ Note that $p$ has degree at most $n_k$. For any $x, y \in \bB_{d_k}$, take $R > 0$ such that $$\max\{|x|,|y|\}<R<1.$$ Harris's range inclusion \cite[Theorem~2]{Harris96} on $p(R\,\cdot)$ gives
    \begin{equation*}
        n_kp(x)+Dp(x)(y-x)\in n_kp(R\overline{\bB_{d_k}}),
        \quad |x|,|y|\le R.
    \end{equation*}

    Let $q_k$ denote the corresponding slice of $Q_k$ and note for a fixed $x \in \bB_{d_k}$ that
    \begin{equation*}
        \cD_x := \{q_k(x) + y \cdot \nabla p(x) \ : \ |y| < 1\} \subset \bC
    \end{equation*}
    is a disk with center $q_k(x)$ and radius $|\nabla p(x)|$. If the radius is $0$, then $\cD_x = \{q_k(x)\}$ is a singleton. In any case, since $P$ is stable, it follows that $0 \notin \cD_x$. Taking $y=0$ gives $q_k(x)\ne0$, and the radius of $\cD_x$ gives $|\nabla p(x)|\le|q_k(x)|$.
    
    By continuity and the fact that the slices and $x$ were arbitrarily chosen, the first bound in \eqref{eq:block.gradient.bound} holds on $\overline\Omega$. Now, the identity $n_kP=Q_k+z^{(k)}\cdot\nabla_k P$ immediately yields the second inequality in \eqref{eq:block.gradient.bound}. Thus, $\eta^{(k)}\cdot\nabla_kP(\eta)=-Q_k(\eta)$ for $\eta\in E(P)$, and so equality in Cauchy--Schwarz gives the final asserted gradient identity.
\end{proof}

By Lemma~\ref{lem:polars}, each active $Q_k$ has smaller total block-degree than $P$ and $E(Q_k)\subseteq E(P)$. The induction will therefore close once $[P]_s$ contains a power of $Q$ using property (4) in Section \ref{sec:properties.Dirichlet-type.space}. The next lemma controls division by $P$ near its boundary zeros. To this end, recall that $E \subseteq \bR^n$ is called \emph{semialgebraic} if it is a finite union of sets of the form
\begin{equation*}
    \{x \in \bR^n \ : \ p_1(x) = \cdots = p_r(x) = 0, \, q_1(x) > 0, \dots, q_s(x) > 0\},
\end{equation*}
where $p_i, q_j$ are polynomials with real coefficients. A function $f : E \to \bR^m$ is called semialgebraic if its graph $\{(x,f(x)) \, : \, x \in E\}$ is semialgebraic.

\begin{lemma}\label{lem:face.error}
    Let $P$ be stable, $E(P)\ne\varnothing$, and $|P|,|Q_k|\le1$
    on $\overline\Omega$. For active $k$, let $\Psi_k=-\nabla_kP/Q_k$,
    and define
    \begin{equation}\label{eq:def.Delta}
        \Delta(z)=\min_{\eta\in E(P)}\max_{k\ \mathrm{active}}
        \max\bigl\{|1-\langle z^{(k)},\eta^{(k)}\rangle|,
                    |1-\Psi_k(z)\cdot\eta^{(k)}|\bigr\}.
    \end{equation}
    For each active $k$, we have $|\Psi_k|\le1$ and $|\cR^j\Psi_k|\lesssim_j |Q|^{-j}$ for $j\in\bN_0$. Also, there is an $\ell \in \bN$ such that
    \begin{equation*}
        |Q(z)|\Delta(z)^\ell\lesssim |P(z)|,\quad z\in\Omega.
    \end{equation*}
\end{lemma}

\begin{proof}
The bound $|\Psi_k|\le1$ follows from \eqref{eq:block.gradient.bound}.
Let $\Omega_A$ be the product of the active blocks and put
$K=\overline\Omega\times\overline{\Omega_A}$. On $K$ define
\begin{equation*}
    f(w,v)=|P(w)|+\max_{k\ \mathrm{active}}
        |1-w^{(k)}\cdot v^{(k)}|
    \qand \cZ=f^{-1}(0).
\end{equation*}
The set $\cZ$ is nonempty, since any $\eta\in E(P)$ gives a point of
$\cZ$ by taking $w=\eta$ and $v^{(k)}=\overline{\eta^{(k)}}$ for each
active $k$.

Fix $(w,v)\in\cZ$. Since $w^{(k)}\cdot v^{(k)}=1$ in every active
block, equality in Cauchy--Schwarz gives $$|w^{(k)}|=|v^{(k)}|=1
\qand v^{(k)}=\overline{w^{(k)}}.$$ Choose $\eta\in\Sigma$ with
$\eta^{(k)}=w^{(k)}$ for every active $k$. Since $P$ depends only on
the active blocks, $P(\eta)=P(w)=0$, so $\eta\in E(P)$. For each
active $k$, note that
\begin{align*}
    |1-\ip{z^{(k)},\eta^{(k)}}|
    &=|\ip{w^{(k)}-z^{(k)},\eta^{(k)}}|
      \le |w^{(k)}-z^{(k)}|,\\
    |1-\Psi_k(z)\cdot\eta^{(k)}|
    &=|(v^{(k)}-\Psi_k(z))\cdot\eta^{(k)}|  \le |v^{(k)}-\Psi_k(z)|.
\end{align*}
Both R.H.S. above are bounded by $|(z,\Psi(z))-(w,v)|$. Therefore, using
this $\eta$ in the minimum defining $\Delta(z)$ and then taking
the infimum over $(w,v)\in\cZ$ gives
\begin{equation*}
    \Delta(z)\le\operatorname{dist}((z,\Psi(z)),\cZ).
\end{equation*}

The compact \L ojasiewicz inequality \cite[Corollary~2.6.7]{BCR98}
states that if $f,g$ are continuous semialgebraic functions on a
compact semialgebraic set and $f^{-1}(0)\subseteq g^{-1}(0)$, then
$|g|^\ell\lesssim |f|$ for some $\ell\in\bN$. Here, take
$g(w,v)=\operatorname{dist}((w,v),\cZ)$ on $K$. Both $f$ and $g$ are
continuous and semialgebraic, and their zero set on $K$ is $\cZ$.
Since $(z,\Psi(z))\in K$, the inequality and the polar identity give
\begin{equation*}
    \Delta(z)^\ell
    \lesssim |P(z)|+\max_{k\ \mathrm{active}}
                  |1-z^{(k)}\cdot\Psi_k(z)|= |P(z)|+\max_{k\ \mathrm{active}}n_k
                  \left|\frac{P(z)}{Q_k(z)}\right|
    \lesssim \frac{|P(z)|}{|Q(z)|}.
\end{equation*}

Finally, successive differentiation of
$Q_k\Psi_k=-\nabla_kP$, using $|Q|\le |Q_k|\le1$, gives
$|\cR^j\Psi_k|\lesssim_j |Q|^{-j}$. Indeed, apply induction on $j$, starting with $|\Psi_k| \leq 1$. The Leibniz rule gives
\begin{equation*}
    Q_k \cR^j \Psi_k = - \cR^j (\nabla_k P) - \sum_{i = 1}^j \binom{j}{i}(\cR^i Q_k)(\cR^{j - i}\Psi_k).
\end{equation*}
Every derivative of $P$ or $Q_k$ is bounded on $\overline{\Omega}$. Since $j - i \leq j -1$ and $0 < |Q| \leq 1$, the induction hypothesis gives
\begin{equation*}
    |\cR^{j-i} \Psi_k| \lesssim_j |Q|^{-(j-i)} \leq |Q|^{-(j-1)}.
\end{equation*}
Dividing by $|Q_k| \geq |Q|$ contributes one more factor of $|Q|^{-1}$, and hence $|\cR^j \Psi_k| \lesssim_j |Q|^{-j}$ holds for all $j$.
\end{proof}

Our final ingredient in the proof of Theorem \ref{mainthm:characterization.of.cyclicity} is the following technical result.

\begin{proposition}\label{prop:product.approximation}
    If $P$ is stable and $E(P)\ne\varnothing$, then
    $Q^{D+1}\in[P]_s$ for all $s\le D-\kappa_*(P)$.
\end{proposition}

\begin{proof}
Scale $P$ so that $|P|,|Q_k|\le1$ on $\overline\Omega$, and let
$\Psi$ and $\Delta$ be as in Lemma~\ref{lem:face.error}. Let
$\Omega_A$ be the product of the active blocks, put
$d=\sum_{k\ \mathrm{active}}d_k$, and let $S_A$ be the
Hardy--Szeg\H{o} kernel on $\Omega_A$. A full tuple in an argument
of $S_A$ is understood to be restricted to the active blocks.
Write $U_v=U_v(E(P))$. At a boundary zero $\eta$ with $Q(\eta) \neq 0$, Lemma \ref{lem:polars} gives $\Psi_k(\eta) = \overline{\eta^{(k)}}$ in every active block, so the following Hardy kernels below are conjugates at $\eta$. For $0<\tau,t\le1/2$, set
\begin{align*}
    K_\tau(z,\omega)
    &=\Bigl[\tau^d S_A\bigl((1-\tau)z,\omega\bigr)\Bigr]
      \Bigl[\tau^d S_A\bigl((1-\tau)\Psi(z),\bar\omega\bigr)\Bigr],\\
    L_t(z)
    &=\int_0^t\!\int_{U_\tau} K_\tau(z,\omega)^{D+2}
       \,d\sigma(\omega)\,\frac{d\tau}{\tau^{d+1}}.
\end{align*}

The following estimates allow us to control the denominators appearing in $K_\tau$:
\begin{equation}\label{eq:mod.asymp}
    \max\left\{1-r,\tfrac12|1-w|\right\}
    \le |1-rw|
    \le 1-r+|1-w|, \quad |w| \leq 1, \, 0 \leq r < 1.
\end{equation}
In particular, for $r=1-\tau$, we have $|1-(1-\tau)w|\asymp \tau+|1-w| \geq \tau$.

\medskip
\noindent\textbf{Step 1. Kernel decay and radial derivatives.}
For $\omega\in U_\tau$, choose $\eta\in E(P)$ with
$\varrho(\omega,\eta)<\sqrt\tau$. Applying the blockwise triangle
inequality through $\omega^{(k)}$ in \eqref{eq:def.Delta} gives
\begin{equation}\label{eq:Delta.bound.proof}
    \Delta(z)\le 2\tau+2\max_{k\ \mathrm{active}}
        \max\bigl\{|1-\ip{z^{(k)},\omega^{(k)}}|,
                    |1-\Psi_k(z)\cdot\omega^{(k)}|\bigr\}.
\end{equation}
Each block factor of $K_\tau$ has the form
$(\tau/(1-(1-\tau)w))^{d_k}$ with $|w|\le1$, using
$|\Psi_k|\le1$ from Lemma~\ref{lem:face.error}. By
\eqref{eq:mod.asymp}, its modulus is at most $1$. For the factor
attaining the maximum in \eqref{eq:Delta.bound.proof}, the
denominator has modulus comparable to $\tau+|1-w|$ and hence
bounded below by a constant times $\tau+\Delta(z)$. Since
$d_k\ge1$, this gives the pointwise bound $|K_\tau(z,\omega)|\lesssim (1+\Delta(z)/\tau)^{-1}$. By \eqref{eqn:Hardy--Szego.kernel}, for $x\in\overline{\Omega_A}$,
\begin{equation*}
    \int_\Sigma |\tau^d S_A((1-\tau)x,\omega)|^2\,d\sigma(\omega)
    =\frac{\tau^{2d}}{\prod_{k\ \mathrm{active}}
      (1-(1-\tau)^2|x^{(k)}|^2)^{d_k}}
    \le\tau^d,
\end{equation*}
since each denominator before exponentiation is at least $\tau$.
Cauchy--Schwarz on the two kernel factors therefore gives the integral bound $\int_{U_\tau}|K_\tau|\,d\sigma\le\tau^d$.

Note that
\begin{equation*}
    \cR K_\tau = (1- \tau) K_\tau \sum_{k \text{ active}} d_k \left( \frac{\ip{z^{(k)}, \omega^{(k)}}}{1 - (1-\tau)\ip{z^{(k)},\omega^{(k)}}} + \frac{(\cR \Psi_k)(z) \cdot \omega^{(k)}}{1 - (1 - \tau) \Psi_k(z) \cdot \omega^{(k)}} \right).
\end{equation*}
Since the denominators have modulus at least $\tau$ and $|\cR \Psi_k| \lesssim |Q|^{-1}$, it follows that $|\cR K_\tau| \lesssim (|Q|\tau)^{-1} |K_\tau|$. By the product and chain rules, \eqref{eq:mod.asymp} and the radial bounds from Lemma~\ref{lem:face.error} give $$|\cR^j(K_\tau^{D+2})|\lesssim (|Q|\tau)^{-j}|K_\tau|^{D+2}, \quad \forall 0 \leq j \leq D.$$
Using the integral bound on one factor $|K_\tau|$ from above, and the pointwise bound on the remaining $D+1$ gives
\begin{equation}\label{eq:paired.kernel.radial.scale}
    \frac1{\tau^d}\int_{U_\tau}
       |\cR^j(K_\tau^{D+2})|\,d\sigma
    \lesssim |Q|^{-j}\tau^{-j}(1+\Delta/\tau)^{-D-1},
    \quad 0\le j\le D.
\end{equation}
On compact subsets of $\Omega$, $|Q|$ and $\Delta$ are bounded
away from zero. Thus the right-hand side is locally uniformly
$O(\tau^{D+1-j})$ near zero and is integrable against $d\tau/\tau$.
It follows that $L_t$ is holomorphic, and differentiation under
the integral gives
\begin{equation}\label{eq:paired.kernel.radial.integral}
    |\cR^j L_t|
    \lesssim |Q|^{-j}\int_0^t
       \tau^{-j}(1+\Delta/\tau)^{-D-1}\,\frac{d\tau}{\tau},
    \quad 0\le j\le D.
\end{equation}

\medskip
\noindent\textbf{Step 2. The lower bound for $\re L_t$.}
Fix $z\in\Omega$ and choose $\eta\in E(P)$ attaining the minimum
in \eqref{eq:def.Delta}. For $x=z^{(k)}$ or
$\overline{\Psi_k(z)}$, we have
$|1-\ip{x,\eta^{(k)}}|\le\Delta(z)$. Splitting
$\omega^{(k)}=\eta^{(k)}+(\omega^{(k)}-\eta^{(k)})$ gives
\begin{align*}
    |\ip{x-\eta^{(k)},\omega^{(k)}}|
    &\le |1-\ip{x,\eta^{(k)}}|
       +|x-\eta^{(k)}|\,|\omega^{(k)}-\eta^{(k)}|\\
    &\le \Delta+
        2\sqrt{\Delta\,|1-\ip{\eta^{(k)},\omega^{(k)}}|},
\end{align*}
where the last inequality uses
$|x-\eta^{(k)}|^2\le2|1-\ip{x,\eta^{(k)}}|$ and the same
inequality for $\omega^{(k)}$ in place of $x$.
By \eqref{eq:mod.asymp},
\begin{align*}
 \frac{(1-\tau)|\ip{x-\eta^{(k)},\omega^{(k)}}|}
       {|1-(1-\tau)\ip{\eta^{(k)},\omega^{(k)}}|} &\lesssim
    \frac{\Delta+2\sqrt{\Delta\,|1-\ip{\eta^{(k)},\omega^{(k)}}|}}
         {\tau+|1-\ip{\eta^{(k)},\omega^{(k)}}|}\\
    &\le \frac{\Delta}{\tau}+\sqrt{\frac{\Delta}{\tau}}.
\end{align*}
Indeed, the first term is bounded using the $\tau$ in the
denominator, while the second uses
$2\sqrt{\tau a}\le\tau+a$ with
$a=|1-\ip{\eta^{(k)},\omega^{(k)}}|$.

For $\tau\ge C\Delta$ with $C$ sufficiently large, each denominator
in $K_\tau(z,\omega)$ is therefore its value at the replacement
$z^{(k)}=\eta^{(k)}$, $\Psi_k(z)=\overline{\eta^{(k)}}$ times
$1+O(\sqrt{\Delta/\tau})$. Taking the fixed negative powers
and multiplying preserves this relative error. The two kernel
factors at the replacement values are conjugates, so
\begin{equation*}
    K_\tau(z,\omega)^{D+2}
    =|\tau^d S_A((1-\tau)\eta,\omega)|^{2(D+2)}
      \bigl(1+O(\sqrt{\Delta/\tau})\bigr).
\end{equation*}
Taking $C$ large enough that the error has modulus at most
$1/2$ gives a lower bound for the real part by half the displayed
power of the kernel.

Since $P$ is independent of the inactive coordinates, those
coordinates of any $\xi\in E(P)$ may be replaced by the
corresponding coordinates of $\omega\in\Sigma$. Hence
\begin{equation}\label{eq:paired.kernel.active.distance}
    \varrho(\omega,E(P))^2
    =\min_{\xi\in E(P)}\max_{k\ \mathrm{active}}
        |1-\ip{\omega^{(k)},\xi^{(k)}}|.
\end{equation}
In particular, the set of $\omega\in\Sigma$ with
$|1-\ip{\eta^{(k)},\omega^{(k)}}|<\tau$ in every active block
is contained in $U_\tau$. This set has measure comparable to
$\tau^d$ by the blockwise Kor\'anyi ball volume estimate.
On this set, each denominator of $S_A((1-\tau)\eta,\omega)$
has modulus at most $2\tau$ by \eqref{eq:mod.asymp}, so
$$\tau^d|S_A((1-\tau)\eta,\omega)|\ge\frac{1}{2^d}.$$ As a result,
\begin{equation*}
    \frac1{\tau^d}\int_{U_\tau}
        \re\bigl(K_\tau(z,\omega)^{D+2}\bigr)\,d\sigma(\omega)
    \ge c_0>0,\quad C\Delta\le\tau\le t.
\end{equation*}

By \eqref{eq:paired.kernel.radial.scale} with $j=0$, the
contribution to $L_t$ from $0<\tau<\min\{t,C\Delta\}$ is
bounded in modulus by a constant times
\begin{equation*}
    \int_0^{C\Delta}(\tau/\Delta)^{D+1}\,\frac{d\tau}{\tau}
    =\frac{C^{D+1}}{D+1}.
\end{equation*}
If $t\ge C\Delta$, the remaining integral has real part at least
$c_0\log(t/(C\Delta))$. If $\Delta\le t<C\Delta$, then
$\log(t/\Delta)\le\log C$. Therefore
\begin{equation}\label{eq:paired.kernel.log.lower}
    \re L_t\ge c_0\log(t/\Delta)-O(1),\quad \Delta\le t.
\end{equation}

\medskip
\noindent\textbf{Step 3. Division by $P$.}
Let $\ell$ be as in Lemma~\ref{lem:face.error}, put
$b=\ell(D+1)$, and choose $c>0$ with $cc_0\ge b$. For
$\Delta\le t$, \eqref{eq:paired.kernel.log.lower} gives
$|e^{-cL_t}|\lesssim(\Delta/t)^b$. For $1\le j\le D$, split
\eqref{eq:paired.kernel.radial.integral} at $\tau=\Delta$ to get
\begin{align*}
    |\cR^j L_t|
    &\lesssim |Q|^{-j}
      \left(\int_0^\Delta \tau^{-j}(\tau/\Delta)^{D+1}
                  \,\frac{d\tau}{\tau}
            +\int_\Delta^t \tau^{-j}\,\frac{d\tau}{\tau}\right)\\
    &\le |Q|^{-j}\Delta^{-j}
      \left(\frac1{D+1-j}+\frac1j\right)
\end{align*}
up to the constant in \eqref{eq:paired.kernel.radial.integral}.
For $\Delta>t$ and $0\le j\le D$, that estimate gives
\begin{equation*}
    |\cR^j L_t|
    \lesssim |Q|^{-j}\int_0^t
        \tau^{-j}(\tau/\Delta)^{D+1}\,\frac{d\tau}{\tau}
    =\frac{(|Q|t)^{-j}}{D+1-j}(t/\Delta)^{D+1}.
\end{equation*}
In particular, $L_t$ is bounded when $\Delta>t$, so
$|e^{-cL_t}|\lesssim1$ and
$$|1-e^{-cL_t}|\le c|L_t|e^{c|L_t|}\lesssim(t/\Delta)^{D+1}.$$

For $j\ge1$, each term in $\cR^j e^{-cL_t}$ is a constant times
\begin{equation*}
    e^{-cL_t}\prod_{a=1}^m \cR^{j_a}L_t,
    \quad j_a\ge1,\quad j_1+\cdots+j_m=j.
\end{equation*}
When $\Delta\le t$, its modulus is at most a constant times
$(\Delta/t)^b(|Q|\Delta)^{-j}$. This is also bounded by a
constant times $(|Q|t)^{-j}$ since $b\ge j$. When $\Delta>t$,
its modulus is at most a constant times
$(|Q|t)^{-j}(t/\Delta)^{m(D+1)}$. Here $m\ge1$, so the
factor $(t/\Delta)^{m(D+1)}$ is at most $(t/\Delta)^{D+1}$.
Also, $t^{-j}(t/\Delta)^{D+1}\le t^{-b}\Delta^{b-j}$.
Combining these bounds with the bounds at order zero gives
\begin{align}
    |\cR^j e^{-cL_t}|
    &\lesssim |Q|^{-j}t^{-b}\Delta^{b-j},
    \label{eq:paired.kernel.exponential}\\
    |\cR^j(1-e^{-cL_t})|
    &\lesssim |Q|^{-j}t^{-j}(1+\Delta/t)^{-D-1},
    \quad 0\le j\le D.
    \label{eq:paired.kernel.exponential.error}
\end{align}
In the second estimate, the derivatives for $j\ge1$ agree up to
sign with those of $e^{-cL_t}$, and for $j=0$ and $\Delta\le t$
we use $$|1-e^{-cL_t}|\le1+|e^{-cL_t}|\lesssim1.$$

Set $H_t=Q^{D+1}e^{-cL_t}/P$. In each term of
$\cR^i(Q^{D+1})$, at least $D+1-i$ factors of $Q$ are
undifferentiated. Since the differentiated factors are bounded
and $|Q|\le1$, this gives
$$|\cR^i(Q^{D+1})|\lesssim |Q|^{D+1-i}, \quad \forall 0 \leq i \leq D.$$
Similarly, applying the Leibniz rule to $P(1/P)=1$ gives
$|\cR^k(1/P)|\lesssim |P|^{-k-1}$ by induction. Indeed, for
$k\ge1$,
\begin{equation*}
    |\cR^k(1/P)|
    \lesssim |P|^{-1}\sum_{a=1}^k |P|^{-(k-a+1)}
    \lesssim |P|^{-k-1},
\end{equation*}
since the radial derivatives of $P$ are bounded and $|P|\le1$.

For $i+j+k\le D$, a Leibniz term with $i,j,k$ derivatives on
$Q^{D+1},e^{-cL_t},1/P$, respectively, is therefore bounded
using \eqref{eq:paired.kernel.exponential} and
$|Q|\Delta^\ell\lesssim |P|$ by a constant times
\begin{equation*}
    |Q|^{D+1-i}|Q|^{-j}t^{-b}\Delta^{b-j}|P|^{-k-1}
    \lesssim t^{-b}|Q|^{D-i-j-k}\Delta^{\ell(D-k)-j}.
\end{equation*}
Both exponents on the right are nonnegative, since $i+j+k\le D$
and $\ell\ge1$. As $|Q|\le1$ and $\Delta\le2$, the radial
derivatives of $H_t$ through order $D$ are bounded for fixed
$t$. Thus $H_t\in\sD_{2D}$ by \eqref{eq:radial.norm.identity}.

\medskip
\noindent\textbf{Step 4. The approximation error.}
For any $0 < t \leq 1/2$, define $$F_t=Q^{D+1}-PH_t=Q^{D+1}(1-e^{-cL_t}).$$
By the Leibniz rule and \eqref{eq:paired.kernel.exponential.error},
for $0\le j\le D$,
\begin{align*}
    |\cR^j F_t|
    &\lesssim (1+\Delta/t)^{-D-1}
       \sum_{i=0}^j |Q|^{D+1-i}|Q|^{-(j-i)}t^{-(j-i)}\\
    &\lesssim t^{-j}(1+\Delta/t)^{-D-1},
\end{align*}
since $|Q|^{D+1-i}|Q|^{-(j-i)}=|Q|^{D+1-j}\le1$ and
$t^{-(j-i)}\le t^{-j}$.

For $\zeta\in\Sigma$ and $0<\rho<1$, the first term in the
maximum in \eqref{eq:def.Delta}, together with
\eqref{eq:mod.asymp} and \eqref{eq:paired.kernel.active.distance},
gives
\begin{align*}
    \Delta(\rho\zeta)
    &\ge \frac12\min_{\eta\in E(P)}\max_{k\ \mathrm{active}}
          |1-\ip{\zeta^{(k)},\eta^{(k)}}|\\
    &=\frac12\varrho(\zeta,E(P))^2.
\end{align*}
Put $\gamma=D-\kappa_*(P)$. By \eqref{eq:sigma.tube.bound},
$\sigma(U_v)\lesssim v^\gamma$ for $0<v\le1/2$, and the bound
extends to $v>1/2$ since $\sigma(\Sigma)=1$.

Integrating the radial estimate and taking the supremum over
$0<\rho<1$ as in \eqref{eq:homogeneous.norm} gives
\begin{equation*}
    \|\cR^j F_t\|_0^2
    \lesssim t^{-2j}\int_\Sigma
       \left(1+\frac{\varrho(\zeta,E(P))^2}{t}\right)^{-2D-2}d\sigma(\zeta).
\end{equation*}
For fixed $\zeta$, we have $\zeta\in U_{tv}$ precisely when
$v>\varrho(\zeta,E(P))^2/t$. Hence
\begin{align*}
    \int_\Sigma \left(1+\frac{\varrho(\zeta,E(P))^2}{t}\right)^{-2D-2}d\sigma(\zeta)
    &=(2D+2)\int_\Sigma\int_{\frac{\varrho(\zeta,E(P))^2}{t}}^\infty
              \frac{dv}{(1+v)^{2D+3}}\,d\sigma(\zeta)\\
    &=(2D+2)\int_0^\infty
              \frac{\sigma(U_{tv})}{(1+v)^{2D+3}}\,dv\\
    &\lesssim t^\gamma\int_0^\infty
              \frac{v^\gamma}{(1+v)^{2D+3}}\,dv
      \\
      &\lesssim t^\gamma.
\end{align*}
The second equality is Tonelli's theorem, and the last integral
converges since $0<\gamma\le D$. This proves
$$\|\cR^j F_t\|_0^2\lesssim t^{\gamma-2j}, \quad \forall 0 \leq j \leq D.$$

In particular, $\|F_t\|_0\lesssim t^{\gamma/2}\to0$. Since
$|F_t(0)|\le\|F_t\|_0$, \eqref{eq:radial.norm.identity} gives
$\|F_t\|_{2D}\lesssim t^{\gamma/2-D}$. By
\eqref{eq:scale.interpolation},
\begin{align*}
    \|F_t\|_\gamma
    &\le \|F_t\|_0^{1-\frac{\gamma}{2D}}
         \|F_t\|_{2D}^{\frac{\gamma}{2D}}\\
    &\lesssim
    t^{\frac{\gamma}{2}(1-\frac{\gamma}{2D})
       +(\frac{\gamma}{2}-D)\frac{\gamma}{2D}}\\
    &=1.
\end{align*}
Hardy norm convergence gives pointwise convergence to zero, so
the weak convergence criterion from property (3) in
Section~\ref{sec:properties.Dirichlet-type.space} gives
$F_t\rightharpoonup0$ in $\sD_\gamma$ as $t\downarrow0$.
Since $H_t\in\sD_{2D}\subseteq\sD_\gamma$, polynomial density
and boundedness of multiplication by $P$ give
$PH_t\in[P]_\gamma$. The weak limit of $PH_t$ is $Q^{D+1}$,
so $Q^{D+1}\in[P]_\gamma$ by weak closedness. Finally,
$\|\cdot\|_s\le\|\cdot\|_\gamma$ gives $[P]_\gamma\subseteq[P]_s$
for every $s\le\gamma$.
\end{proof}

\begin{proof}[\textbf{Proof of Theorem~\ref{mainthm:characterization.of.cyclicity}}]

Assume henceforth that $P$ is stable on $\Omega$. That a cyclic polynomial must be stable is already established in property (4) of Section \ref{sec:properties.Dirichlet-type.space}. If $E(P)=\varnothing$, then
$\inf_{\Sigma}|P|>0$. Uniform convergence of $P(\rho\,\cdot)$ to $P$ on
$\Sigma$ and the blockwise maximum principle on $\rho\Omega$ give a
uniform bound for $1/P$ as $\rho\uparrow1$. Thus $P$ has no zeros on
$\overline\Omega$, and $1/P\in\Hol(\overline\Omega)\subseteq\sD_s$.
It follows that $1\in[P]_s$ for every $s$.

Suppose $E(P)\ne\varnothing$ and put $\gamma=D-\kappa_*(P)$.
For $s\le\gamma$, induct on the block-total degree, starting with
nonzero constants. By Lemma~\ref{lem:polars}, each active $Q_k$ is
stable, has smaller block-total degree, and satisfies
$E(Q_k)\subseteq E(P)$. The definition of the measure exponent gives
$\kappa_*(Q_k)\le\kappa_*(P)$, including the empty-set case by
convention. Hence each active $Q_k$ is cyclic in $\sD_s$ by induction,
and so is $Q^{D+1}$ by property (4) in Section~\ref{sec:properties.Dirichlet-type.space}.
Proposition~\ref{prop:product.approximation} gives $Q^{D+1}\in[P]_s$,
so $[P]_s=\sD_s$.

For the converse, take $s>\gamma$ and let $\mu$ be the probability
measure from Theorem~\ref{mainthm:covering-measure.estimates}.
Let $f=\sum_{n}f_n$ be the homogeneous expansion of a polynomial and, for $t = 2^{-j-1}$, set
\[
    f^{[j]}=\sum_{2^j\le n+1<2^{j+1}}f_n \qand g_j=\sum_{2^j\le n+1<2^{j+1}}(1-t)^{-n}f_n.
\]
The factors $(1-t)^{-n}$ are uniformly bounded because $nt\le1$
and $t\le1/2$. By \eqref{eq:homogeneous.norm},
$\|g_j\|_0\lesssim\|f^{[j]}\|_0$. Since
$g_j((1-t)\eta)=f^{[j]}(\eta)$, the Hardy-trace estimate \eqref{eq:Hardy-trace} gives
\[
    \|f^{[j]}\|_{L^2(\mu)}
    \lesssim t^{-\gamma/2}\|g_j\|_0
    \lesssim 2^{j\gamma/2}\|f^{[j]}\|_0.
\]
Using that $\mu$ is a probability measure and applying
Cauchy--Schwarz to the sum over $j$, we obtain
\begin{equation*}
    \left|\int_{E(P)}f\,d\mu\right|
    \lesssim \sum_{j\ge0}2^{j\gamma/2}\|f^{[j]}\|_0\le
    \left(\sum_{j\ge0}2^{j(\gamma-s)}\right)^{1/2}
    \left(\sum_{j\ge0}2^{js}\|f^{[j]}\|_0^2\right)^{1/2}
    \lesssim_s \|f\|_s.
\end{equation*}
Thus $f\mapsto\int_{E(P)}f\,d\mu$ extends continuously to $\sD_s$, vanishes on polynomial multiples of $P$, and takes the value
$1$ on the constant function $1$. Hence $1\notin[P]_s$.
\end{proof}

Corollaries~\ref{cor:ball.characterization} and~\ref{cor:polydisk.characterization}
follow at once from the following result. Recall that for a polynomial $P$, $m_P := \dim_\bR E(P)$ with the convention that $m_P = -\infty$ if $E(P) = \varnothing$.

\begin{lemma}\label{lem:classical.covering}
For a stable $P$, $\kappa_*(P) = m_P/2$ on $\bB_D$ and $\kappa_*(P) = m_P$ on $\bD^D$.
\end{lemma}

\begin{proof}
    Suppose $E(P)\ne\varnothing$ and let $F_j : \cQ_j \to \Sigma$ be as in Section \ref{sec:proof.Theorem.A}. By the inverse function theorem, each $F_j(\cQ_j)$ has a projection onto $m_j$ real coordinates with nonempty interior. Every projection onto more than $m_j$ coordinates has measure zero by Sard's theorem (see \cite[Theorem~6.10]{Lee13}). Since $E(P)=\bigcup_j F_j(\cQ_j)$ is a finite union, $m_P=\max_j m_j$.

    On the ball, fix $x\in\cQ_j$ and $y\in\bR^{m_j}$. The analytic arc $a(t)=F_j(x+ty)$ extends holomorphically near $0$, and $P\circ a=0$ there by the identity theorem. By stability, $|a(it)|^2\ge |a(0)|^2=1$ for real $t$ near $0$. Differentiating at this minimum gives
    \[
        0=\left.\frac{d}{dt}|a(it)|^2\right|_{t=0} =-2\im\ip{DF_j(x)y,F_j(x)}.
    \]
    Thus $A_j(x)=0$, and the rank formula gives $\kappa_*(P)=m_P/2$.

    On the polydisk, differentiating $|F_j^{(k)}|^2=1$ in each scalar block gives
    \[
        DF_j^{(k)}(x)y=iF_j^{(k)}(x)(A_j(x)y)_k.
    \]
    Hence $\ker A_j(x)=\ker DF_j(x)=\{0\}$, so $h_j=m_j$ and $\kappa_*(P)=m_P$.
\end{proof}

\section{Examples}\label{sec:special.case.example}

\begin{example} In this example, we show that the real dimension of the distinguished boundary zero set does
not determine the cyclicity endpoint on a general polyball. The
following two polynomials on $\bB_2\times\bD$ both have circles
as their distinguished boundary zero sets. Let
    \[
        P(z_1,z_2,u)=2-z_1^2-z_2^2-u
        \qand Q(z_1,z_2,u)=1-z_1u.
    \]
    On $\bB_2\times\bD$ we have
    $|z_1^2+z_2^2+u|\leq |z_1|^2+|z_2|^2+|u|<2$ and
    $|z_1u|<1$, so both polynomials are stable. At a
    distinguished boundary zero of $P$, the two terms
    $z_1^2+z_2^2$ and $u$ have modulus at most one and sum
    to $2$, so both equal $1$. Since
    $\re(z_j^2)=|z_j|^2-2(\im z_j)^2$ and
    $|z_1|^2+|z_2|^2=1$, this forces $z_1,z_2$ to be real.
    At a distinguished boundary zero of $Q$, the equation
    $z_1u=1$ gives $|z_1|=1$, $z_2=0$, and $u=\bar z_1$.
    Conversely, these conditions give zeros, so
    \begin{align*}
        E(P)&=\{\eta_\theta=(\cos\theta,\sin\theta,1):
                      \theta\in\bR\},\\
        E(Q)&=\{\xi_\theta=(e^{i\theta},0,e^{-i\theta}):
                      \theta\in\bR\}.
    \end{align*}
    In the first parametrization the block inner products
    between the points at $\theta$ and $\phi$ are
    $\cos(\theta-\phi)$ and $1$. In the second they are
    $e^{i(\theta-\phi)}$ and $e^{-i(\theta-\phi)}$. Thus
    \begin{align*}
        \varrho(\eta_\theta,\eta_\phi)^2
            &=1-\cos(\theta-\phi),\\
        \varrho(\xi_\theta,\xi_\phi)^2
            &=|1-e^{i(\theta-\phi)}|.
    \end{align*}
    For $|\alpha|\leq\pi$ we have
    $|1-e^{i\alpha}|=2|\sin(\alpha/2)|\asymp|\alpha|$ and
    $1-\cos\alpha=\tfrac12|1-e^{i\alpha}|^2\asymp\alpha^2$.
    Therefore, on $E(P)$ the condition
    $\varrho(\eta_\theta,\eta_\phi)^2<t$ cuts out an interval
    of angles of length comparable to $\sqrt t$ modulo $2\pi$.
    On $E(Q)$ the corresponding interval has length comparable
    to $t$. The measures induced by $d\theta/(2\pi)$ give the
    same powers for the masses of these balls.

    To see that these powers are optimal, let $E$ denote either
    circle and $\nu$ its probability measure, with
    $\nu(B_t(\eta))\asymp t^a$ for $\eta\in E$. If a ball
    $B_t(\xi)$ centered anywhere on $\Sigma$ meets $E$ at $\eta$,
    then $B_t(\xi)\subseteq B_{4t}(\eta)$ by the triangle inequality.
    Hence $\nu(B_t(\xi))\lesssim t^a$ for every center $\xi$.
    If another probability measure $\mu$ supported on $E$
    satisfies \eqref{eq:polynomial.measure} with exponent
    $\kappa$, then symmetry of the balls and Tonelli's theorem give
    \[
        t^a\lesssim
        \int_E\nu(B_t(\eta))\,d\mu(\eta)
        =\int_E\mu(B_t(\xi))\,d\nu(\xi)
        \lesssim t^\kappa .
    \]
    Letting $t\to0$ gives $\kappa\leq a$.

    We conclude that $m_P=m_Q=1$, while $\kappa_*(P)=1/2$
    and $\kappa_*(Q)=1$. Since $D=3$,
    Theorem~\ref{mainthm:characterization.of.cyclicity} gives
    cyclicity exactly for $s\leq5/2$ for $P$ and $s\leq2$ for $Q$.
\end{example}

\begin{example} Here, we observe that enlarging the domain can raise or lower the cyclicity endpoint. Both possibilities occur for $Q=1-z_1u$ on $\bB_2 \times \bD$ from the preceding example. First replace the disk block by a ball and put
    \[
        \widetilde Q(z_1,z_2,u_1,u_2)=1-z_1u_1
        \qquad\text{on }\bB_2\times\bB_2.
    \]
    This polynomial is stable since $|z_1u_1|<1$ in the
    interior. At a distinguished boundary zero,
    $z_1u_1=1$ forces $|z_1|=|u_1|=1$, hence $z_2=u_2=0$.
    Consequently,
    \[
        E(\widetilde Q)
        =\{(e^{i\theta},0,e^{-i\theta},0):\theta\in\bR\}.
    \]
    The two block inner products are again
    $e^{i(\theta-\phi)}$ and $e^{-i(\theta-\phi)}$.
    Thus the measure computation for $Q$ also gives
    $\kappa_*(\widetilde Q)=1$. The new complex dimension is
    $4$, so the endpoint increases from $2$ to $4-1=3$.

    Next consider $Q$ on $\bD^3$, which strictly contains
    $\bB_2\times\bD$. It remains stable since $|z_1u|<1$
    on $\bD^3$. On the distinguished boundary, $z_1u=1$
    gives $u=\bar z_1$, while $z_2$ ranges freely over $\bT$.
    Hence
    \[
        E(Q)=\{(e^{i\theta},e^{i\phi},e^{-i\theta}):
                                      \theta,\phi\in\bR\}.
    \]
    Projection onto the first two complex coordinates identifies
    this set with $\bT^2$, so its real dimension is $2$.
    Corollary~\ref{cor:polydisk.characterization} gives the
    endpoint $3-2=1$ on $\bD^3$.
\end{example}

\begin{example} The block decomposition of a polyball can also affect cyclicity even when the polynomial, the ambient dimension, and the distinguished boundary zero set are all fixed. To this end, let $P(z)=2-(z_1+z_2)(z_3+z_4)$ on $\bC^4$.
    Consider the two polyballs with blocks
    $(z_1,z_2),(z_3,z_4)$ and $(z_1,z_3),(z_2,z_4)$,
    respectively.

    For the first decomposition, Cauchy--Schwarz gives
    \[
        |(z_1+z_2)(z_3+z_4)|
        \leq 2|(z_1,z_2)|\,|(z_3,z_4)|\leq 2
    \]
    on the closure, with strict inequality in the interior.
    At a zero, equality requires both blocks to have norm one
    and each block to be a multiple of $(1,1)$.
    Thus $z_1=z_2$, $z_3=z_4$, and all four coordinates have
    modulus $1/\sqrt2$.

    For the second decomposition, the estimate on the closure is
    \[
        |(z_1+z_2)(z_3+z_4)|
        \leq \tfrac12|(z_1,z_3)+(z_2,z_4)|^2\leq 2.
    \]
    Again the last inequality is strict in the interior.
    Equality holds there on the closure only if
    $(z_1,z_3)=(z_2,z_4)$ and both vectors have norm one.
    Equality in the first inequality gives
    $|z_1+z_2|=|z_3+z_4|$, so $|z_1|=|z_3|=1/\sqrt2$.

    Thus $P$ is stable on both polyballs. Under either
    decomposition, a zero on the closure satisfies
    $4z_1z_3=2$ and hence $z_3=\bar z_1$. Conversely, each
    point of the circle
    \[
        \{\eta_\theta=(a,a,\bar a,\bar a):
                       a=e^{i\theta}/\sqrt2,\ \theta\in\bR\}
    \]
    belongs to both distinguished boundaries and is a zero of $P$.

    For the first decomposition, the two block inner products
    at $\eta_\theta,\eta_\phi$ are $e^{i(\theta-\phi)}$ and
    $e^{-i(\theta-\phi)}$. For the second, both inner products
    equal
    $\tfrac12(e^{i(\theta-\phi)}+e^{-i(\theta-\phi)})
     =\cos(\theta-\phi)$.
    Thus the squared distances are $|1-e^{i(\theta-\phi)}|$
    and $1-\cos(\theta-\phi)$, respectively. The circle
    computations in the first example give exponents $1$
    and $1/2$. The corresponding cyclicity endpoints are
    $4-1=3$ and $4-1/2=7/2$.
\end{example}

\begin{example} We end our discussion with an explicit endpoint for a family of polynomials over any polyball.

    Let $1\leq h_k\leq d_k$ for $1\leq k\leq r$, and put
    \[
        F_k(z^{(k)})=\sum_{j=1}^{h_k}(z_j^{(k)})^2,
        \qquad P(z)=1-\prod_{k=1}^r F_k(z^{(k)}).
    \]
    Since $|F_k(z^{(k)})|\leq|z^{(k)}|^2<1$ in each ball,
    $P$ is stable. At a distinguished boundary zero of $P$,
    the product equals one and each $|F_k|\leq1$, so every
    $|F_k|=1$. Choose $\theta_k$ such that
    $F_k(z^{(k)})=e^{2i\theta_k}$. Then
    \[
        1=\sum_{j=1}^{h_k}
               \re\bigl((e^{-i\theta_k}z_j^{(k)})^2\bigr)
          \leq\sum_{j=1}^{h_k}|z_j^{(k)}|^2
          \leq |z^{(k)}|^2=1.
    \]
    Equality in the first inequality makes each
    $e^{-i\theta_k}z_j^{(k)}$ real for $j\leq h_k$.
    Equality in the second makes all coordinates with
    $j>h_k$ vanish. It follows that $E(P)$ consists exactly
    of the points
    \[
        z^{(k)}=e^{i\theta_k}(x_k,0),\qquad
        x_k\in\mathbb S_{\mathbb R}^{h_k-1},\qquad
        \sum_{k=1}^r\theta_k=0\pmod{\pi},
    \]
    where $\mathbb S_{\mathbb R}^{h_k-1}$ is the real unit sphere
    in $\bR^{h_k}$. Indeed, at such a point
    $F_k(z^{(k)})=e^{2i\theta_k}$, so the angle condition makes
    their product equal to one.

    Locally, one can choose the argument of each
    $F_k(z^{(k)})$ continuously, which fixes $\theta_k$ and
    then determines $x_k$. Thus local coordinates consist
    of $r-1$ independent angles and local coordinates on the
    real spheres. This gives real dimension
    $r-1+\sum_k(h_k-1)$. When $h_k=1$, the corresponding
    real sphere consists of two points and contributes zero
    to the local dimension.

    A variation $(\dot\theta_k,\dot x_k)$ in these coordinates
    gives
    \[
        \dot z^{(k)}
        =e^{i\theta_k}\bigl(i\dot\theta_k(x_k,0)+(\dot x_k,0)\bigr).
    \]
    Because $|x_k|=1$ and $\dot x_k\cdot x_k$ is real,
    \[
        \im\ip{\dot z^{(k)},z^{(k)}}
        =\im\bigl(i\dot\theta_k+\dot x_k\cdot x_k\bigr)
        =\dot\theta_k.
    \]
    The possible angle variations are exactly the tuples
    satisfying $\sum_k\dot\theta_k=0$. Thus the map in
    \eqref{eq:intrinsic.exponent.formula} has rank $r-1$.
    Substituting the dimension and rank gives
    \[
        \kappa_*(P)
        =\frac12\left(r-1+\sum_{k=1}^r(h_k-1)+r-1\right)
        =r-1+\frac12\sum_{k=1}^r(h_k-1).
    \]
    By Theorem~\ref{mainthm:characterization.of.cyclicity},
    $P$ is cyclic in $\sD_s(\Omega)$ exactly for
    $s\leq D-r+1-\tfrac12\sum_{k=1}^r(h_k-1)$.
\end{example}

\end{document}